\documentclass[12pt]{article}
 \usepackage{amsfonts}  
 \usepackage{graphicx}
\title{Regularity and uniqueness of the first eigenfunction for singular fully non linear operators. }

\author{ I. Birindelli F. Demengel}

\date{}

\catcode`@=11 \@addtoreset{equation}{section} \catcode`@=12

\newtheorem{theo}{Theorem}[section]
\newtheorem{prop}[theo]{Proposition}
\newtheorem{rema}[theo]{Remark}
\newtheorem{defi}[theo]{Definition}
\newtheorem{cor}[theo]{Corollary}
\newtheorem{lemme}[theo]{Lemma}
\newtheorem{exe}[theo]{Example}

\def\R{{\rm I}\!{\rm  R}}
\def\grad{\nabla}
\def\ml{{\mathcal L}}
\def\Fl{{\cal F}}

\begin{document}
\maketitle
\section{Introduction}
The concept of principal eigenvalue for boundary value problems of elliptic operators has been extended, in the last decades, to quasi-nonlinear and fully-nonlinear equations (somehow "abusing" the name of eigenvalue) see \cite{A,SA,OT,QS,IY,BD3} etc.. In all the cases we know, two features of the operators are requested {\em homogeneity} and {\em ellipticity}. 

The meta-definition of these principal eigenvalues could be the following: Given a zero order  and odd
operator $H$ with the same homogeneity than the second order elliptic operator $F$, and given a domain $\Omega$, 

$\lambda$ is a principal  eigenvalue if there exists a non trivial solution of constant sign of the problem
$$\left\{\begin{array}{lc}
F(x,\grad u, D^2u)+\lambda H(u)=0 &\mbox{in}\quad\Omega\\
u=0 &\mbox{on}\quad\partial \Omega;
\end{array}
\right.
$$
not surprisingly, that function will be called the "eigenfunction".

\bigskip
Of course when $F$ is a second order linear elliptic operator and $H(u)=u$, the principal eigenvalue
is just the first eigenvalue in the "classical" sense, it is well known that it is both {\em simple} and {\em isolated}. 

It is interesting to notice that when the operator is not "odd" with respect to the Hessian, it can be expected that there are two principal eigenvalues, one corresponding to a positive eigenfunction and one corresponding to a negative eigenfunction, and this is the case for example when $F$ is one of the Pucci operators e.g. for some $0<a<A$
$$F(D^2u):={\cal M}_{a,A}^+(D^2u):=a\sum_{e_i<0}e_i+A\sum_{e_i>0}e_i$$
where $e_i$ denote the eigenvalues of the matrix $D^2u$  (see \cite{BEQ}). In this cases, or more in general  when the operator is uniformly elliptic and homogenous of degree 1, the principal eigenvalues have been proved to be simple and isolated (see \cite{QS,IY,P}).

When the operator is "quasi linear" but in divergence form , for example in the case of the $p$-Laplacian, it is well known that the principal eigenvalue can be defined through the Rayleigh quotient and it was proved, independently  by  Anane \cite{A} and by Otani and Teshima \cite{OT}, that it is simple and isolated (see also \cite{LP2}). The variational structure plays a key role there. On the other hand, for the $\infty$-Laplacian (see \cite{J}), the question of the simplicity of the eigenvalue is still open.

\bigskip

The cases treated in this paper concern operators that have the "homogeneity" of the $p$-Laplacian, but are "fully-nonlinear" and hence are not variational. In previous works, we proved the existence of the principal eigenvalues for this large class of operators and many features related to them (\cite{BD1},$\dots$,\cite{BD5}). The inspiration for these definitions and results was the acclaimed work of Berestycki, Nirenberg and Varadhan \cite{BNV} where the eigenvalue for linear elliptic operators in general bounded domains was defined through the maximum principle. 

The main questions left open  in our previous works were: are these eigenvalues "simple" ? are they "isolated"?

We shall now proceed to describe the results obtained in this note  but for the sake of comprehension we shall do it for an operator that  exemplifies well the cases treated here (the general conditions and hypothesis will be given in the next section).
For some $0\geq \alpha>-1$, and some H\"older continuous function $h$ of exponent $<1+\alpha$, let:
$$F[u]:=F(x,\grad u, D^2u)=|\grad u|^\alpha{\cal M}_{a,A}^+(D^2u) +h(x)\cdot\grad u |\grad u|^\alpha.$$
Suppose that $\Omega$ is a bounded, smooth domain of $\R^N$. Then we can define
$$\lambda^+:=\sup\{\lambda; \exists \phi>0,\quad F[\phi]+\lambda\phi^{1+\alpha}\leq 0\quad\mbox{in }\Omega\},$$
$$\lambda^-:=\sup\{\lambda; \exists \phi<0,\quad F[\phi]+\lambda\phi|\phi|^{\alpha}\geq 0\quad\mbox{in }\Omega\}.$$
It is clear that the inequalities are meant in the "viscosity sense" adapted to these non smooth operators (see the next section for a precise definition).

In a much more general context and for any $\alpha>-1$ we proved that these are well defined;
that for any $\lambda<\lambda^+$ the maximum principle holds i.e. considering the problem
\begin{equation}\label{000}
\left\{
\begin{array}{lc}
F[u]+\lambda u|u|^{\alpha}=0 &\mbox{in}\quad\Omega\\
u=0 &\mbox{on}\quad\partial \Omega,
\end{array}
\right.
\end{equation}
if $u$ is a viscosity subsolution  of (\ref{000}) then $u\leq 0$ in $\Omega$ and for any $\lambda\leq \min\{\lambda^+,\lambda^-\}$ and any continuous $f$ there exists a solution of
\begin{equation}\label{001}
\left\{
\begin{array}{lc}
F[u]+\lambda u|u|^{\alpha}=f &\mbox{in}\quad\Omega\\
u=0 &\mbox{on}\quad\partial \Omega,
\end{array}
\right.
\end{equation}
which is Lipschitz continuous.

Of course we also proved that there exists $\phi^+$ and $\phi^-$ respectively positive and negative eigenfunctions in the sense that e.g. for $\lambda=\lambda^+$ there exists $\phi^+>0$ viscosity solution of (\ref{000}).

One of the question we raise here is , if $\phi>0$ is another solution of (\ref{000}) with $\lambda=\lambda^+$, is it true that there exists $t>0$ such that $\phi^+=t\phi$? The answer is yes for any domain $\Omega$ such that $\partial\Omega$ has only one connected component.
When $\partial\Omega$ has two connected components we can prove the result when  $N=2$.
We also recall that for any $\alpha>-1$ the simplicity of the eigenvalue has been proved in 
\cite{BD5} in the case of radial solutions for any $N$.

It is clear that these results are somehow equivalent to a "strong comparison principle"
i.e. it is equivalent to know that if two solutions are one above the other in some open set ${\cal O}$ and they "touch" at some point of ${\cal O}$ then they coincide in ${\cal O}$.
Now in Proposition \ref{prop1}, we prove such a result  when at least one of the solutions has the gradient away from zero in ${\cal O}$.

This restriction implies that  in order to apply Proposition \ref{prop1} we need to know that there is 
some subset ${\cal O}$ of $\Omega$ where this condition is satisfied. Naturally the Hopf's lemma together with a $C^1$ regularity  is the right ingredient since it guarantees that the gradient is away from zero in a neighborhood of $\partial\Omega$.
This explains why we start by proving a $C^{1,\beta}$ regularity result, which is interesting in itself but 
furthermore is essential in the proof of the simplicity of the principal eigenvalue.

This is done through a fixed point theorem.
Let us mention that in general the tools to prove regularity are the Alexandrov-Bakelman-Pucci inequality and some "sub-linearity" of the operator. In a recent paper Davila, Felmer and Quaas have proved ABP \cite{DFQ} for singular fullynonlinear operators, but in this case it does not seem to be useful to prove regularity because the difference of a sub and super solution may not be a sub-solution of an elliptic equation.

The dimensional restriction is due to the fact that when there are two  connected components of the boundary of $\Omega$, we use Sard's Theorem. A famous counterexample of Whitney
shows that Sard's theorem doesn't hold if the functions are only $C^1$. It seems that the least regularity that can be asked is $W^{N,p}$ (see e.g. \cite{Fi}), and since the solutions are in $W^{2,p}$ we have to take $N=2$.

Other important results concerning eigenvalues are given as a consequence of simplicity. In particular we prove that there are no eigenfunctions that change sign for $\lambda=\lambda^{\pm}$. Further results include  the strict monotonicity of the eigenvalue with respect to the inclusion of domains. And finally that the eigenvalues are isolated.

The paper is organized as follows: In the next section we give the precise hypothesis concerning the operator and we recall the known results concerning singular operators. In section three we prove the $C^{1,\beta}$ regularity of the solutions.
Section four is devoted to the strict comparison principle and simplicity of the principal eigenvalues. We end the paper with other properties of the eigenvalues.

\section{ Assumptions and known results}

Let $\Omega$ be a bounded regular set of $\R^N$. 

 Let us recall what we mean by {\it viscosity solutions}, adapted to
our context.

     \begin{defi}\label{def1}

 Let $\Omega$ be a bounded domain in
$\R^N$, let $g$ be a continuous function on $\Omega\times \R$, then
$v$,   continuous on $\overline{\Omega}$ is called a viscosity super-solution (respectively sub-solution) 
of
$F(x,\grad u,D^2u)=g(x,u)$ if for all $x_0\in \Omega$, 

-Either there exists an open ball $B(x_0,\delta)$, $\delta>0$  in $\Omega$
on which 
$v= cte= c
$ and 
$0\leq g(x,c)$, for all $x\in B(x_0,\delta)$
(respectively 
$0\geq g(x, c)$ for all $x\in B(x_0,\delta)$) 

-Or
 $\forall \varphi\in {\mathcal C}^2(\Omega)$, such that
$v-\varphi$ has a local minimum (respectively local maximum) at  $x_0$ and $\grad\varphi(x_0)\neq
0$, one has
$$
F( x_0,\grad\varphi(x_0),
 D^2\varphi(x_0))\leq g(x_0,v(x_0)).
$$
(respectively  
$$F( x_0,\grad\varphi(x_0),
 D^2\varphi(x_0))\geq g(x_0,v(x_0)).$$

A	 viscosity solution is a function which is both a super-solution and a sub-solution.
\end{defi}

\begin{rema}
 When $F$ is continuous in $p$,  and $F(x,0,0)=0$, this  definition is equivalent to the classical definition
of viscosity solutions, as in the User's guide \cite{CIL}.
\end{rema}
We now 
state the assumptions satisfied by the operator $F$. Let $S$ be the set of $N\times N$ symmetric matrices, and let  $\alpha \in ]-1,0[$. Then $F$ defined on $\Omega \times \R^N\setminus\{0\}\times S $ satisfes 
\begin{equation}\label{deff}
F(x, p, M) = |p|^\alpha  (\tilde F(x, M)+h(x)\cdot p).
\end{equation}
on $\tilde F$ we suppose

\begin{itemize}
\item[(F)] $\tilde F(x,tM)=t\tilde F(x,M)$ for any $t\in \R^+$, for any $M\in S$,  and any $N\in S$ such that $N\geq 0$ there exist $A\geq a>0$
\begin{equation}\label{eqaAF}
a tr(N)\leq \tilde F(x,M+N)-\tilde F(x,M) \leq A tr(N).
\end{equation}
Furthermore $(x, M)\mapsto \tilde F(x,M)$ is continuous.

\item[(J)]
 There exists a
continuous function $  \omega$ with $\omega (0) = 0$, such that if
$(X,Y)\in S^2$ and 
$\zeta\in \R^+$ satisfy
$$-\zeta \left(\begin{array}{cc} I&0\\
0&I
\end{array}
\right)\leq \left(\begin{array}{cc}
X&0\\
0&Y
\end{array}\right)\leq 4\zeta \left( \begin{array}{cc}
I&-I\\
-I&I\end{array}\right)$$
and $I$ is the identity matrix in $\R^N$,
then for all  $(x,y)\in \R^N$, $x\neq y$
$$\tilde F(x, X)-\tilde F(y,  -Y)\leq \omega
(\zeta|x-y|^2).$$
\end{itemize}
On $h$ we suppose that :

(H)  $h$ is H\"older continuous of exponent $1+\alpha$.

\begin{rema}
Since $\alpha <0$,  if $\tilde F$ satisfies condition (J), 
then so does $F$ i.e. for $X$, $Y$ and $\zeta$ as above 
$$ F(x, \zeta (x-y), X)- F(y, \zeta(x-y),  -Y)\leq \omega
(\zeta|x-y|^2).$$
\end{rema}

   \begin{exe} We suppose that $h$ satisfies (H).
   
   1) Let  $0< a < A$ and ${\cal M}_{a,A}^+(M)$ be the Pucci's operator 
   ${\cal M}_{a,A}^+ (M) = Atr(M^+)-a tr(M^-)$ where $M^\pm$ are the positive and negative part    of $M$, and 
    ${\cal M}_{a,A}^-(M)=-  {\cal M}_{a,A}^+ (-M)$. 
       Then $F$ defined as  
   
   $$F(x, p,M) = |p|^\alpha( {\cal M}_{a,A}^\pm (M) + h(x) \cdot p)$$
 satisfies the assumptions.
   
    2) More generally let $B(x)$ be some matrix with Lipschitz coefficients which is invertible for 
all $x\in \Omega$. Let us consider $A(x) = B^\star B (x)$ and the operator 
    
    $F(x, p,M) = |p|^\alpha(tr(A(x) (M)) + h(x) \cdot p)$,
Then  $\tilde F$ satisfies  (F) and (J),  arguing as  in \cite{BD1}, example 2.4.

    \end{exe}
 We begin to recall some of the results obtained in \cite{BD3} which will be needed in this article.

 \begin{theo}\label{thcomp1}
 Suppose that  $F$,  $h$  are as above and $c$  is continuous and bounded and satisfies  $c\leq 0$. 

Suppose that $f_1$ and $f_2$ are continuous and bounded and that 
$u$ and
$v$  satisfy
\begin{eqnarray*}
 F(x, \nabla u, D^2 u)+c(x)|u|^\alpha u& \geq & f_1\quad 
\mbox{in}\quad \Omega \\ 
F(x,  \grad v,D^2 v)
+ c(x) |v|^{\alpha}v & \leq & f_2 \quad  \mbox{in}\quad
\Omega  \\ 
u \leq  v &&   \quad  \mbox{on}\quad \partial\Omega.
\end{eqnarray*}
 Suppose that $f_2< f_1$, then $u \leq v$ in $\Omega$. 
Moreover if $c<0$ in $\Omega$,  and $f_2\leq f_1$ the result still holds. 

\end{theo}

This comparison theorem  allows to prove,  using the existence of sub- and supersolutions constructed with the aid of the distance function to $\partial \Omega$, together with Perron's method adapted to our context,  the following existence's result : 

 \begin{theo}\label{exi}
  Suppose that  $F$,  $h$ and $c$ are as above and that $c\leq 0$. 
  Suppose that $f$ is continuous and bounded, then there exists a  continuous solution to 
  
  \begin{equation}\label{j09}
\left\{ \begin{array}{lc}
  F(x, \nabla u, D^2 u)  +c(x)|u|^\alpha u= f\quad& \  
\mbox{in}\quad \Omega \\
 u = 0 & {\rm on} \ \partial \Omega
  \end{array}\right.
\end{equation}
  If furthermore  $f\leq 0$, $u \geq 0$, and if $f\geq 0$, $u\leq 0$. 
  Moreover if $c<0$ in $\Omega$, the solution is unique. 
  \end{theo}
  
   \begin{rema}\label{remhopf}
Recall that by the  Hopf principle: if ${\cal O}$ is a smooth bounded domain, and
 $u$ is  a solution of 
    \begin{equation}\label{123}
  F(x, \nabla u, D^2 u) \leq 0 
\end{equation}
in ${\cal O}$, such that $u>c$ inside ${\cal O}$ and $u(\bar x)=c$ at some boundary point of ${\cal O}$ then
$"{\partial u\over \partial n}(\bar x)"<0$. 

In particular this implies that a non constant solution of (\ref{123}) in a domain $\Omega$ has no interior minimum.
    \end{rema}

   We also recall some  regularity results 
\begin{prop}\label{prophold}

Suppose that $\tilde F$ satisfies (F), (J) and $h$ is continuous.  
Let $f$ be some continuous  function in $\overline{\Omega}$.
 Let $u$ be a viscosity non-negative bounded solution of 
\begin{equation}\label{eq4.1}
\left\{
\begin{array}{lc}
F(x, \nabla u, D^2u)=f & \ {\rm in}\
\Omega\\
u=0 &  \ {\rm on}\ \partial\Omega.
\end{array}
\right.
\end{equation}
Then,  for
any
$\gamma\leq 1$ there exists some constant
$C$ which depends only on $|f|_\infty$, $|h|_\infty$ and $|u|_\infty$ such that :

$$|u(x)-u(y)|\leq C|x-y|^\gamma$$
for any
$(x,y)\in\overline\Omega^2$.
\end{prop}

 This Proposition implies  some compactness result for sequences of solutions  which will be used in  section 3. 
 
 We shall also need in the proof of Theorem \ref{uni} the following comparison principle, proved in  
 \cite{BD4}.
 
 \begin{theo}\label{complambda}
 
Suppose that $\tilde F$ satisfies (F), (J), and  $h$ satisfies (H),   
that  $c$ is  continuous and
bounded.   
 Suppose that $0\leq \lambda < \bar\lambda (\Omega) $, that $c+\lambda$  is positive on $\overline{\Omega}$ and that $u$ and $v$ are respectively  positive continuous super and sub solutions of

$$|\nabla u|^\alpha \left(\tilde F(x,D^2 u) + h(x).\nabla u\right) +(c(x)+\lambda)u^{1+\alpha}=  0$$

1) If $u\geq v>0$ on $\partial\Omega$ then $u\geq v$ in $\Omega$.

\noindent 2) If $u> v$ on $\partial \Omega$ then $u> v$ on $\overline\Omega$.

\end{theo}

 \section{Regularity}
 In this section for $\alpha \leq 0$ we establish that the solutions of (\ref{eq4.1}) are $C^{1,\beta}$; this will be a  consequence of the known regularity results in the case $\alpha = 0$. 
 
 We expect that the $C^{1,\beta}$ regularity of the solutions is true for more general operators i.e.
 operators that are only homogeneous and singularly elliptic but not necessarily of the form given in (\ref{deff}). As an example, in the second subsection we illustrate this for a class of operators which does not satisfy
 (\ref{deff}) and is somehow close to the $\infty$-Laplacian, though not as degenerate.
 
  \subsection{Regularity result for $\alpha \leq 0$}
 
 To prove the regularity result announced (which will be stated precisely in Corollary \ref{rglr}) , we  remark  that the solutions of the Dirichlet problem can be obtained  as a fixed point of some operator acting in  ${\cal C}^{1, \beta}$ for all $\beta < \infty$. 

\bigskip
  We define 
 ${\mathcal L}(\Omega):=\{u\in {\cal C}_o(\overline\Omega);\ u \mbox{ is Lipschitz}\ \}$ with the norm $\displaystyle |u|_{\ml}:=\sup_{x,y\in\Omega, x\neq y} \frac{|u(x)-u(y)|}{|x-y|} $. 

 And we suppose here that $g:\R^{+}\rightarrow(\R^+)^\star$ is continuous and decreasing, while $l\mapsto lg(|l|)$ is increasing on $\R$.  

 \begin{prop}\label{propT}
Suppose that    $\tilde F$, $h$  satisfy respectively  (F), (J)  and (H). 
 Let $f\in L^\infty(\Omega)$,  let  $T_\epsilon$  be the operator $\ml\mapsto \ml$ such that $T_\epsilon u=v $,  where $v$ is the unique solution of 
  $$ \left\{ \begin{array}{lc}
 \tilde F(x,D^2 v) + h(x)\cdot \nabla v =  (f+\epsilon g(|u|) u ){1\over g(|\nabla u|)} & {\rm in} \ \Omega \\
   v = 0 & {\rm on} \ \partial \Omega
   \end{array}\right.$$
   then there exists $\epsilon_o$ small enough in order that  for $\epsilon < \epsilon_o$,  
   $T_\epsilon $ has a non trivial fixed point in $\ml$.
     \end{prop}
\begin{rema}
The solution is taken in the sense of $L^p$ viscosity solutions,see \cite{CCKS} and \cite{NW}.
\end{rema}

Before giving the proof of the Proposition \ref{propT} we shall prove the main result of this section i.e. the following 

\begin{cor}\label{rglr} Suppose that $\tilde F$, $h$  satisfy respectively  (F)  (J) and (H) and that
$f\in {\cal C} (\overline{\Omega})$. Let $u$ be a solution  in the  sense of Definition \ref{def1} of  
  $$ \left\{ \begin{array}{lc}
 |\nabla u|^\alpha \left( \tilde F(x,D^2 u) + h(x)\cdot \nabla u \right)=  f & {\rm in} \ \Omega \\
u = 0 & {\rm on} \ \partial \Omega
\end{array}\right.$$
Then  $u\in {\cal C}^{1, \beta}(\overline{\Omega})$ for all $\beta \in(0,1)$.
\end{cor}
 
\begin{rema}

 When the operator $M\mapsto F(x, M)$ is convex or concave, the solutions are ${\cal C}^{2,\beta}$. 
 \end{rema}

{\em Proof of Corollary \ref{rglr}.} 
 We prove  first the announced regularity results for the equation 
 $$ \left\{ \begin{array}{lc}
 |\nabla u|^\alpha \left( \tilde F (x,D^2 u) + h(x)\cdot \nabla u \right)-\epsilon |u|^\alpha u =  f & {\rm in} \ \Omega \\
u = 0 & {\rm on} \ \partial \Omega
\end{array}\right.$$
 for $\epsilon$ small enough.

Let $\delta>0$ and suppose that $v_\delta$ is a fixed point for $T_\epsilon$ with 
$g(l)=(|l|^2+\delta^2)^{\frac{\alpha}{2}}$ which exists by Proposition \ref{propT}. We shall prove that it is the unique solution, in the sense of Definition \ref{def1}, of
\begin{equation}\label{7mai} \left\{ \begin{array}{lc}
 (|\nabla u|^2+\delta^2)^\frac{\alpha}{2} \left( \tilde F (x,D^2 u) + h(x)\cdot \nabla u \right)-\epsilon (\delta^2+ |u|^2)^{\alpha \over 2} u =  f & {\rm in} \ \Omega \\
u = 0 & {\rm on} \ \partial \Omega.
\end{array}\right.
\end{equation}
Observe first  that since $\tilde F$ satisfies (F), by standard regularity results 
  $v_\delta\in {\cal C}^{1,\beta}(\overline\Omega)$ for any $\beta\in (0,1)$ (see Evans \cite{E}, Caffarelli- Cabr\'e \cite{CC} and Winter \cite{NW}). 

We prove only the super-solution case. Let $\bar x\in \Omega$ and $\varphi$  be in ${\cal C}^2$ such that $(v_\delta-\varphi)(x) \geq (v_\delta-\varphi)(\bar x)=0$, with $\nabla \varphi(\bar x)\neq 0$ . Since $v_\delta$ is ${\cal C}^{1,\beta}$,  $\nabla \varphi (\bar x) =\nabla v_\delta(\bar x)$ and since $v_\delta $ is a sub-solution of the fixed point equation: 
\begin{eqnarray*}
 \tilde F (\bar x,D^2 
 \varphi(\bar x))+ h(\bar x)\cdot \nabla \varphi (\bar x)& \geq&  ( f(\bar x)+\epsilon ( |v_\delta(\bar x)|^2+ \delta^2)^{\alpha\over 2}  v_\delta(\bar x)) (|\nabla v_\delta|^2(\bar x)+\delta^2)^{-\frac{\alpha}{2}}\\
 &
= & ( f+\epsilon ( |\varphi(\bar x)|^2+ \delta^2)^{\alpha\over 2}  \varphi(\bar x))( |\nabla \varphi|^2(\bar x)+\delta^2)^{-\alpha\over 2}
\end{eqnarray*}
i.e.
 
 $$(|\nabla \varphi |^2(\bar x)+ \delta^2)^{\alpha\over 2} 
  (\bar x) \left( \tilde F (\bar x,D^2 \varphi(\bar x)) + h(\bar x)\cdot \nabla \varphi(\bar x) \right)-\epsilon( |\varphi (\bar x)|^2+\delta^2)^{\alpha\over 2}  (\bar x) \varphi (\bar x) \geq f.$$
 and then $v_\delta$ is a supersolution of (\ref{7mai}). 
 
We suppose now that for some $c_o$,   $v_\delta(x)=c_o$ for $x$ in a neighborhood of $\bar x$. From the equation we get  that 
$-\epsilon (\delta^2+c_o^2)^{\alpha\over 2}  c_o\leq f$  in a neighborhood of $\bar x$ which is the same condition required in the Definition \ref{def1}.

Proceeding similarly for the subsolution we have obtained that  $v_\delta$ is a solution of (\ref{7mai}) and there exists $C$ depending only on the structural constants of $\tilde F$ and $h$ such that
\begin{eqnarray*}
\|v_\delta\|_{C^{1,\beta}}&\leq &C(|f|_\infty+\epsilon |v_\delta |_\infty ^{1+\alpha})(|\nabla v_\delta|_\infty^2+\delta^2)^{-\frac{\alpha}{2}}\\
&\leq& C(|f|_\infty+\epsilon|v_\delta |_\infty^{1+\alpha})(|\nabla v_\delta|_\infty^2+1)^{-\frac{\alpha}{2}}.
\end{eqnarray*}

Arguing as in (\cite{BD3}), one can prove that  $v_\delta$ is uniformly Lipschitz; this implies both that the bounds
on the $C^{1,\beta}$ norm don't depend on $\delta$ and that we can pass to the limit in the equation (\ref{7mai}). Stability of viscosity solutions implies that $v_\delta$ converges to
the unique  solution of

$$ \left\{ \begin{array}{lc}
 |\nabla u|^\alpha \left( \tilde F(x,D^2 u) + h(x)\cdot \nabla u \right)-\epsilon |u|^\alpha u=  f & {\rm in} \ \Omega \\
u = 0 & {\rm on} \ \partial \Omega.
\end{array}\right.$$
 We now consider the case  $\epsilon = 0$. We write the equation in the following way
   $$
     |\nabla u|^\alpha     \left(  \tilde F (x,D^2 u) + h(x)\cdot \nabla u \right)-\epsilon |u|^\alpha u
      = f-\epsilon |u|^\alpha u.
  $$
Replacing $f$ by $f-\epsilon |u|^\alpha u$ since $u \in L^\infty$,  it gives the $C^{1,\beta}$ regularity  of $u$ and ends the proof.
 
 \begin{rema} Under the hypothesis of Corollary \ref{rglr} if $g:\Omega\times\R\rightarrow \R$ is continuous and bounded, then any solution of 
 
 $$ \left\{ \begin{array}{lc}
 |\nabla u|^\alpha \left( \tilde F(x,D^2 u) + h(x)\cdot \nabla u \right)+ g(x,u)=  f & {\rm in} \ \Omega \\
u = 0 & {\rm on} \ \partial \Omega
\end{array}\right.$$
is in ${\cal C}^{1, \beta}(\overline{\Omega})$ for any $\beta\in (0,1)$.
\end{rema}
  {\em  Proof of Proposition \ref{propT}.}
Clearly $\ml(\Omega)$, equipped with the norm $|.|_\ml$ is a Banach space.
 It is well known that for any $k\in L^\infty$ there exists a unique solution $v$  of 
         
$$ \left\{ \begin{array}{lc}
 \tilde F (x,D^2 v) + h(x)\cdot \nabla v =  k& {\rm in} \ \Omega \\
   v = 0 & {\rm on} \ \partial \Omega.
\end{array}\right.$$
Furthermore $v$ is Lipschitz continuous and there exists some constant $c$  which depends only on the structural data such that 
   $$|v|_\ml \leq c |k|_\infty.$$
   Let $\epsilon$ be small enough in order that $\epsilon c\  d_\Omega<1$ where $d_\Omega$ is the diameter of $\Omega$.

 Let  $G$ be  the reciprocal function of $l\mapsto lg(l)$. 
Suppose first  that $ d_\Omega\leq 1$. Then we choose
 $f_o={c|f|_\infty\over 1-\epsilon c }$ and
   $${\cal B} = \left\{ u \in \ml( \Omega), \mbox{such that} \ |u|_{\ml}\leq G\left(f_o\right)\right\}. $$
${\cal B}$ is a convex and compact set of  $\ml$.  We use the fact that $l\mapsto lg(|l|)$ and $\frac{1}{g}$ are increasing to obtain that,
if $u\in {\cal B}$, 
\begin{eqnarray*}
|T_\epsilon(u)|_\ml &\leq &c \left(|f|_\infty + \epsilon \ g(d_\Omega G(f_o))d_\Omega G(f_o) \right) {1\over g\left(G(f_o)\right) }\\
&\leq  &c \left(|f|_\infty + \epsilon \ f_o \right) \frac{G(f_o)}{f_o}=G(f_o).
\end{eqnarray*}
On the other hand if $ d_ \Omega\geq 1$, then we choose
 $f_o={c|f|_\infty\over 1-\epsilon d_\Omega c }$ and
   ${\cal B}$ as before.
  
Hence if $u\in {\cal B}$, using that $g$ is decreasing:
\begin{eqnarray*}
|T_\epsilon(u)|_\ml &\leq&c \left(|f|_\infty + \epsilon \  g(d_\Omega G(f_o))d_\Omega G(f_o) \right) {1\over g\left(G(f_o)\right)} \\
&\leq &c \left(|f|_\infty + \epsilon \ d_\Omega\ g(G(f_o)) G(f_o) \right) {1\over g\left(G(f_o)\right)} \\
&\leq&c \left(|f|_\infty + \epsilon d_\Omega \ f_o \right) \frac{G(f_o)}{f_o}=G(f_o). 
\end{eqnarray*}
In both cases  $T(B)\subset B$ and furthermore $T$ is continuous and compact,  then Schauder's fixed point theorem implies the result.

     \subsection{Other operators}
  
 We now present an example to which the results of the previous section  can be extended to but  which does not satisfy the previous assumptions.
      
     Let  $q\geq 0$, and  $\alpha \leq 0$, 
     $${\cal F}(p, M) = |p|^\alpha \left ( tr M+ q \langle {Mp\over |p|}, {p\over |p|}\rangle\right) := |p|^\alpha \tilde F(M,p).$$
\begin{prop} Suppose that $u$ is a viscosity solution of 
\begin{equation}\label{plap}\left\{\begin{array}{lc}
\Fl(\grad u,D^2u)=f & \mbox{in}\ \Omega\\
u=0& \mbox{on}\ \partial\Omega
\end{array}\right.
\end{equation}
with $f\in {\cal C}(\overline\Omega)$. Then $u$ is in $C^{1,\beta}(\overline\Omega)$ for any $\beta\in (0,1)$.
\end{prop}     
{\em Proof}  

Recall that the $(q+2)$-Laplacian is defined by
$$\Delta_{q+2}u=|\grad u|^q\left(\Delta u +q\left(D^2u(\frac{\grad u}{|\grad u|}),\frac{\grad u}{|\grad u|}\right)\right).$$

We first prove a fixed point property.    For $\epsilon>0$ and $\delta>0$ let $T_{lap}$ be the map  
     $u\mapsto v$ where 
      $v$ is the solution of 
      $$\left\{ \begin{array}{lc}
      \Delta_{q+2} (v)= (f+\epsilon (|u|^2+\delta^2)^{\alpha\over 2}  u)(|\nabla u|^2+\delta^2)^{\frac{q-\alpha}{2}}&\ {\rm in} \ \Omega\\
       u=0 & {\rm on} \ \partial \Omega.
       \end{array}\right.$$
Recall that by regularity results of the $q+2$-Laplacian, there exists $C_{lap}$ such that
\begin{equation}\label{eqTlap}
|T_{lap} (u)|_{{ \cal C}^{1, \beta}} \leq C_{lap} (|f|_\infty+\epsilon(|u|_\infty ^2+\delta^2)^{\alpha\over 2} | u|_\infty )^\frac{1}{q+1} (|\nabla u|_\infty^{2}+\delta^2)^{q-\alpha \over 2(q+1)}.
\end{equation}
For $\epsilon$ small enough, let  $l_\epsilon$  be  a  solution of
$$l_\epsilon^{q+1}=(C_{lap} (|f|_\infty+\epsilon (d_\Omega l_\epsilon)^{\alpha+1})(l_\epsilon^2+\delta^2)^\frac{q-\alpha}{2}, $$
 and define  the  closed convex compact set in $\ml(\Omega)$
         $${\cal B} = \{ u\in \ml(\Omega) , \ |u|_{\ml} \leq l_\epsilon\}.$$
 Then $T_{lap}B\subset B$ and using the Schauder fixed point theorem one gets that there exists $u_\delta \in {\cal C}^{1, \beta}$ which solves 
$$\left\{ \begin{array}{lc}
\Delta_{q+2} (u)= (f+\epsilon(|u|^2+\delta^2)^{\alpha\over 2}  u)(|\nabla u|^2+\delta^2)^{\frac{q-\alpha}{2}}&\ {\rm in} \ \Omega\\
  u=0 & {\rm on} \ \partial \Omega.
       \end{array}\right.$$
As before, we  can show that in fact $u_\delta$ is a solution in the sense of definition \ref{def1}
of     
$$\left\{ \begin{array}{lc}
(|\nabla u|^2+\delta^2)^{\frac{\alpha-q}{2}}\Delta_{q+2} (u)-\epsilon(|u|^2+\delta^2)^{\alpha\over 2}  u= f&\ {\rm in} \ \Omega\\
  u=0 & {\rm on} \ \partial \Omega.
       \end{array}\right.$$
 Arguing as in the proof of H\"older's and Lipschitz's regularity in \cite{BD3}, one can prove some uniform Lipschitz estimates on $v_\delta$. Then  it is clear form  the estimates  in (\ref{eqTlap})  that  the $C^{1,\beta}$ norm of   $v_\delta$ is independent of $\delta$, we can then  let $\delta$ go to zero and obtain that $u_\delta$ converges to a solution $u$  of 
           $$\left\{ \begin{array}{lc}
  |\nabla u|^{\alpha}   \left(\Delta u + q (D^2u(\frac{\grad u}{|\grad u|}),\frac{\grad u}{|\grad u|}) \right)-\epsilon |u|^\alpha u = f&\ {\rm in} \ \Omega\\
       u=0 & {\rm on} \ \partial \Omega.
       \end{array}\right.$$
         We  then obtain the regularity when $\epsilon = 0$ by  writing the equation
                $$\left\{ \begin{array}{lc}
|\grad u|^\alpha\left(     \Delta u + q (D^2u(\frac{\grad u}{|\grad u|}),\frac{\grad u}{|\grad u|})\right)  = f&\ {\rm in} \ \Omega\\
       u=0 & {\rm on} \ \partial \Omega
       \end{array}\right.$$
        when $f\in L^\infty$, under the form 
                    $$\left\{ \begin{array}{lc}
|\grad u|^\alpha\left(     \Delta u + q (D^2u(\frac{\grad u}{|\grad u|}),\frac{\grad u}{|\grad u|}) \right)-\epsilon |u|^\alpha u= f-\epsilon |u|^\alpha u &\ {\rm in} \ \Omega\\
 u=0 & {\rm on} \ \partial \Omega
       \end{array}\right.$$
       Using the fact that $u$ is in $L^\infty$ one gets the regularity result.

  \section{Strict comparison principle and uniqueness of the eigenfunction}
We now assume that $c(x)$ is some continuous bounded function, and we consider the principal  eigenvalues $\lambda^+$, $\lambda^-$  for the operator  
$$|\nabla u|^\alpha\left( \tilde F(x, D^2 u) + h(x)\cdot \nabla u\right)+ c(x) |u|^\alpha u$$ as defined in the introduction. 
We suppose that  conditions (F) and (J) are satisfied by $\tilde F$,  and $h$ satisfies (H). 
 
We prove  the uniqueness result  in the case $\lambda^+$, the changes to bring for ${\lambda^- }$ being obvious.

\begin{theo}\label{uni}

  Suppose that  $\Omega$ is a  bounded regular   domain such that   either $\partial \Omega $ is connected  or $N=2$  and $\partial \Omega$ has at most two connected components. Suppose that $c$ satisfies $c(x)+ \lambda^+> 0$ in $\Omega$.  
  
 If $\psi$ and $\varphi$ are  two positive eigenfunctions for the eigenvalue $\lambda^+$ there exists $t>0$ such that $\psi = t\varphi$. 
 \end{theo}
 
 \begin{rema} The condition $N=2$ is used only in order to apply Sard's theorem. Indeed
 in \cite{Fi}, Sard's theorem is proved to hold for functions  that are in $W^{N,p}(\R^N)$ but we only know that $u$ is in $W^{2,p}$ .
 \end{rema}
 In order to prove Theorem \ref{uni} we shall need a few results concerning comparison principle and applications of  Hopf's principle.
 We begin by  a strong comparison principle inside $\Omega$ for  sub and supersolutions  $u$ and $v$ which coincide on one point  $\bar x$ inside $\Omega$,  where $|\nabla  u|(\bar x)\neq 0$, or $|\nabla  v|(\bar x)\neq 0$:

 \begin{prop}\label{prop1}
 
 Suppose that 
 $u$ and $v$ are respectively nonnegative $C^{1,\beta}$ solutions of

 \begin{eqnarray*}
  F(x,\nabla u,D^2 u)&\leq& f \\
 F(x,\nabla v,D^2 v)&\geq & g  \end{eqnarray*}
 with $f\leq g$.
 Suppose that  ${\cal O}$ is an open connected subset of $\Omega$, such that

1) $u\geq v$, in ${\cal O}$,

2) $\exists\ x_o\in {\cal O}$ such that $u(x_o)> v(x_o)$,

3)  either $|\grad u(x)|\neq 0$
or $|\grad v(x)|\neq 0$ in ${\cal O}$,

\noindent Then $u> v$ in ${\cal O}$.

\end{prop}

{\em Proof of Proposition \ref{prop1}.} 
Using the connectedness of ${\cal O}$ it is sufficient to prove the
result in some ball  containing $x_o$ instead  of ${\cal O}$.

  Suppose by contradiction that there exists some point $x_1$  close to $x_o$, such that $u(x_1) = v(x_1)$ and in the ball $B(x_o, R)$, for  $R = |x_1-x_o|$,  $x_1$ is  the only point in the closure of that ball  on which $u$ and $v$ coincide.     
  One can also assume that $B(x_o, {3R\over 2}) \subset \Omega$.  
  
  We shall prove that there exists some constant  $c>0$ and  some $\delta >0  $ such that $$ u \geq v + \delta  ( e^{-c|x-x_o|}- e^{-3cR\over 2})\equiv v+ w$$  
in the annulus ${R\over 2}\leq |x-x_o| = r\leq {3R\over 2}$. This will  contradict  the fact that $u = v$ on $x_1$.  

 Let $\delta =\displaystyle \min_{|x-x_o|= {R\over 2} }(u-v)$,   so that 
 $$u\geq v+w\quad  \mbox{on} \quad
 \partial\left(B(x_o, {3R\over 2})\setminus {B(x_o, {R\over 2}})\right).$$

Let $\varphi$ be some test function for $v$ from above, a simple calculation on $w$ implies that, if $c \geq {1\over a}({2(2A(N-1) \over R} +2 |h|_\infty)$ then
                
     \begin{eqnarray*} 
     |\nabla \varphi+ \nabla w|^\alpha &\cdot&\left(\tilde F (x, D^2 \varphi+ D^2 w)+h(x)\cdot (\nabla \varphi+ \nabla w)\right) \\
     &\geq&    |\nabla \varphi+ \nabla w|^\alpha \left(\tilde F (x,D^2 \varphi)+h(x)\cdot \nabla \varphi\right)
     \\
     &&+   |\nabla \varphi+ \nabla w|^\alpha\left( {\cal M}^- (D^2 w)-h(x)\cdot \nabla w\right)\\
       &\geq& |\nabla \varphi+ \nabla w|^\alpha { g\over |\nabla \varphi|^\alpha} +\\
       && + |\nabla \varphi+ \nabla w|^\alpha\left(ac^2-Ac({N-1\over r} )-|h|_\infty c \right) \delta e^{-cr}\\
        &\geq& |\nabla \varphi+ \nabla w|^\alpha { g\over |\nabla \varphi|^\alpha} +\\
       && + |\nabla \varphi+ \nabla w|^\alpha\frac{ac^2}{2}  \delta e^{-cr}.
\end{eqnarray*}

We also impose that $\delta < {R L_1e\over 16 }$, which implies in particular that $|\nabla w|\leq {|\nabla \varphi]\over 8}$

    We now use the inequalities 

    $$||\nabla \varphi+\nabla w|^\alpha-|\nabla \varphi|^\alpha |\leq |\alpha | |\nabla w|| \nabla \varphi|^{\alpha-1}\left( {1\over 2}\right)^{\alpha-1}\leq {|\nabla \varphi|^\alpha \over 2}$$
to get  
     \begin{eqnarray*} 
     |\nabla \varphi+ \nabla w|^\alpha \left(\tilde F(x,D^2 \varphi+ D^2 w)\right.&+&\left.h(x)\cdot (\nabla \varphi+ \nabla w)\right) \\
&\geq &  g -| g|_\infty|\grad\varphi|^{-1}|\alpha| 2^{1-\alpha}c\delta e^{-cr}   +L_2^\alpha {ac^2\over 4} \delta e^{-cr}.
\end{eqnarray*}
It is now enough to choose  
$$c=\sup \left\{{2({2A(N-1) \over R} + |h|_\infty) \over a},{2^{4-\alpha} | g|_\infty \over a L_1 L_2^{\alpha} } \right\}$$ 
to finally obtain
\begin{eqnarray*} 
     |\nabla \varphi+ \nabla w|^\alpha \left(\tilde F(x,D^2 \varphi+ D^2 w)\right.&+&\left.h(x)\cdot (\nabla \varphi+ \nabla w)\right) \\
&\geq &   f+{a  c^2\delta L_2^\alpha e^{-cr}\over 8}.
\end{eqnarray*}
i.e.
  $$F(x, \nabla (v+ w), D^2 (v+ w)) > F(x,\nabla u,D^2 u).$$
We  now are in a position to use the comparison principle and  get that 
   $$ u\geq v+w $$
  in   the annulus $B(x_o,\frac{3R}{2})\setminus  B(x_o,\frac{R}{2})$, the desired contradiction.
This ends the proof of Proposition \ref{prop1}.

\bigskip
Another way of formulating this proposition is the following
\begin{cor}\label{234} Suppose that $u$ and $v$ are as in Proposition \ref{prop1} but instead of condition 2) we have that there exists $\bar x\in {\cal O}$ such that $u(\bar x)=v(\bar x)$,  then $u\equiv v$ in ${\cal O}$.
\end{cor}

 Another consequence of Proposition \ref{prop1} is the following :   Let $\partial_{ \nu} u$ be the normal derivative $\nabla u\cdot \vec \nu$ where $\vec \nu $ is the unit outer normal to $\partial \Omega$.  Then 

\begin{prop}\label{chopf}
Suppose that $\Omega$ is a smooth bounded domain of $\R^N$.
Suppose  that $u$ and $v$ are respectively nonnegative ${\mathcal C}^{1,\beta}$ solutions of

 \begin{eqnarray*}
 && F(x,\nabla u,D^2 u)\leq f \\
&& F(x,\nabla v,D^2 v)\geq g \end{eqnarray*}
 with $f\leq g$, and $u = v = 0$ on $\partial \Omega$

1) $u\geq v$ in $\Omega $,

2) there exists $\bar x\in \partial\Omega$ such that $\partial_\nu u(\bar x)=  \partial_\nu v(\bar x)$ 

\noindent then there exists $\epsilon>0$ such that 
$$u\equiv v\ \mbox{in}\ \Omega\setminus\overline\Omega_\epsilon$$
where $\Omega_\epsilon$ is the set of points of $\Omega$ whose distance to  the  connected component of the boundary  which contains $\bar x$ is greater than $\epsilon$.
\end{prop}
{\em Proof.} First with Hopf Principle $|\nabla u| >0$ and $|\nabla v|>0$  on the boundary and by the regularity results  there exist $\delta>0$, $L_1$ and $L_2$ such that $L_1\leq |\nabla u|, |\nabla v|\leq L_2$ in $\Omega\setminus\overline\Omega_\epsilon$. 

If there exists a point $x_1$ of $\Omega\setminus\overline\Omega_\epsilon$ such that 
$u(x_1)=v(x_1)$ by the previous Corollary we have nothing to prove. So we can suppose by contradiction that there exists a ball $B\subset\Omega$ which is tangent to $\partial\Omega$ in $\bar x$ where $u>v$. Let $x_2$ be the center of that ball. Let $R = |\bar x-x_2|$. 

Let  $w= \delta (e^{-c|x- x_2|}-e^{-c|x_2-\bar x|})$  where $\delta$ is chosen such that  $\delta \leq \inf_{|x-x_2|= {R\over 2}} (u-v)$.

Reasoning as in Proposition \ref{prop1} we get that there exists  $c$ such that  
$$u\geq v+  w,$$
in the annulus $B( x_2, R)\setminus B( x_2, {R\over 2})$

 This implies in particular that 
 
  $$|\partial_\nu u(\bar x)|\geq|(\partial_\nu v+ \partial_\nu w)(\bar x)|.$$
Since $|(\partial_\nu v+ \partial_\nu w)(\bar x)|>|\partial_\nu v(\bar x)|$ this leads to
a contradiction. This ends the proof of Proposition \ref{chopf}. 

\bigskip

In the sequel we shall need the following well known result : 

 \begin{lemme}\label{sph}
 Suppose that $O$ is an open  bounded set. There exists some point  on $\partial O$  where $\partial O$ 
 satisfies the interior sphere  condition.  
                   \end{lemme}
                   See \cite{Bo} for more complete results on that property. 
                   
\noindent {\em Proof of Theorem \ref{uni}.}

We suppose first that $\partial \Omega$ is connected. 
 Let $d(x)$ denote the distance to the boundary of $\Omega$.
Suppose that $\psi$ and $\varphi$ are two positive eigenfunctions and
let $\Gamma=\sup\frac{\psi}{\varphi}$. 
This  extremum is  well defined because, using Hopf lemma and the comparison principle, there exist $c_1$ and $c_2$ such that in a neighborhood of $\partial \Omega$:
$$c_1d(x)\leq \psi(x),\varphi (x)\leq c_2 d(x),$$
(see \cite{BD3} for the details).
Moreover this supremum is  achieved on the boundary in the sense that  e.g. there exists a sequence $(x_n)_n$ which goes to the   boundary such that ${\psi\over \varphi} (x_n)\rightarrow \Gamma$: 

Indeed, suppose not, then, there would exist an open set $\Omega^\prime $, $\Omega^\prime \subset \subset \Omega$  such that  on  $\Omega \setminus \overline{\Omega}^\prime $ 
$${\psi\over \varphi} \leq \Gamma-\epsilon$$
for some $\epsilon>0$.  Using the comparison Theorem \ref{complambda} in $\Omega^\prime $, and remarking that $\overline{\lambda} (\Omega) < \overline{\lambda}(\Omega^\prime )$  one would get that 
    $\psi\leq (\Gamma-\epsilon) \varphi$ on $\Omega^\prime $ and finally on all $\Omega$, a contradiction. 

We are in a position to apply
  Corollary \ref{234} and Proposition \ref{chopf}   and obtain that $\psi = \Gamma \phi$ in a neighborhood of $\partial \Omega$. 
  Considering the infimum of the ratio one gets  also that if $ \gamma = \inf \frac{\psi}{ \varphi}$  then $\frac{\psi}{\varphi} = \gamma$ on a neighborhood of the boundary, therefore $\gamma=\Gamma$  and the conclusion follows. 
   
   \bigskip

  We now assume that $N=2$ and    since  $\partial \Omega$ has at most two connected components there exist $\Omega_1$ and $\Omega_2$  simply connected  and smooth such that $\Omega_1\subset \subset \Omega_2$  and $\Omega = \Omega_2 \setminus \overline{\Omega_1}$. 
  
   We  define as before $\Gamma=\sup\frac{\psi}{\varphi}$ and $\gamma=\inf\frac{\psi}{\varphi}$, and arguing as previously there exists some sequence $x_n$ and $i \in \{ 1,2\}$  such that $x_n\rightarrow \partial \Omega_i$ and ${\psi\over \varphi}(x_n)\rightarrow \Gamma$. 
By Corollary \ref{234} and Proposition \ref{chopf} with $f := -(c(x)+ \overline{\lambda }) \psi^{1+\alpha} \leq -(c(x)+ \overline{\lambda })  \varphi^{1+\alpha}:=g$,
$$\psi\equiv\Gamma\varphi\ \mbox{in a neighborhood of}\ \partial \Omega _i.$$
In the same manner we prove that the infimum $\gamma=\inf\frac{\psi}{\varphi}$ is achieved on the other connected component, and also that $\psi = \gamma \varphi$ in a neighborhood of this  part of the frontier. 
 
 We have obtained   that 
$$ \psi\equiv\Gamma_1\varphi\ \mbox{in a neighborhood of}\ \partial \Omega_1$$
and 
$$ \psi\equiv\Gamma_2\varphi\ \mbox{in a neighborhood of}\ \partial \Omega_2.$$
with $\{\Gamma_1, \Gamma_2\} = \{ \Gamma, \gamma\}$.

For $i = 1,2$, let 
$$A_i \mbox{ be the  connected component of   } \{ x, \psi (x) = \Gamma_i \varphi(x),\  \nabla \psi (x)\neq 0\ {\rm or }\  \nabla \varphi (x)\neq 0\}$$  
 whose boundary contains $\partial \Omega_i$.  
  
 $A_i$ is open. 
    Indeed if $x_1\in A_1$ then $\psi(x_1)=\Gamma_1 \varphi(x_1)$ and there is $N_{x_1}$ a neighborhood of $x_1$ where either $\grad \psi\neq 0$ or $\grad \varphi\neq 0$ and then using Corollary \ref{234}, in $N_{x_1}$, $\psi=\Gamma_1 \varphi$ and then $N_{x_1}\subset A_1$.

Observe that if  $\partial A_1\cap \partial \Omega_2\neq \emptyset$ or $\partial A_2\cap \partial \Omega_1\neq \emptyset$, this ends the proof of Theorem \ref{uni} since it would imply that $\Gamma=\gamma$. Hence we suppose that these intersections are empty.

    Let $G_i := \overline{A_i}$;   $K_i := \partial G_i \cap \Omega$,  and  $M_i := \Omega \setminus G_i $. 
 \begin{center}  \includegraphics[scale=0.4]{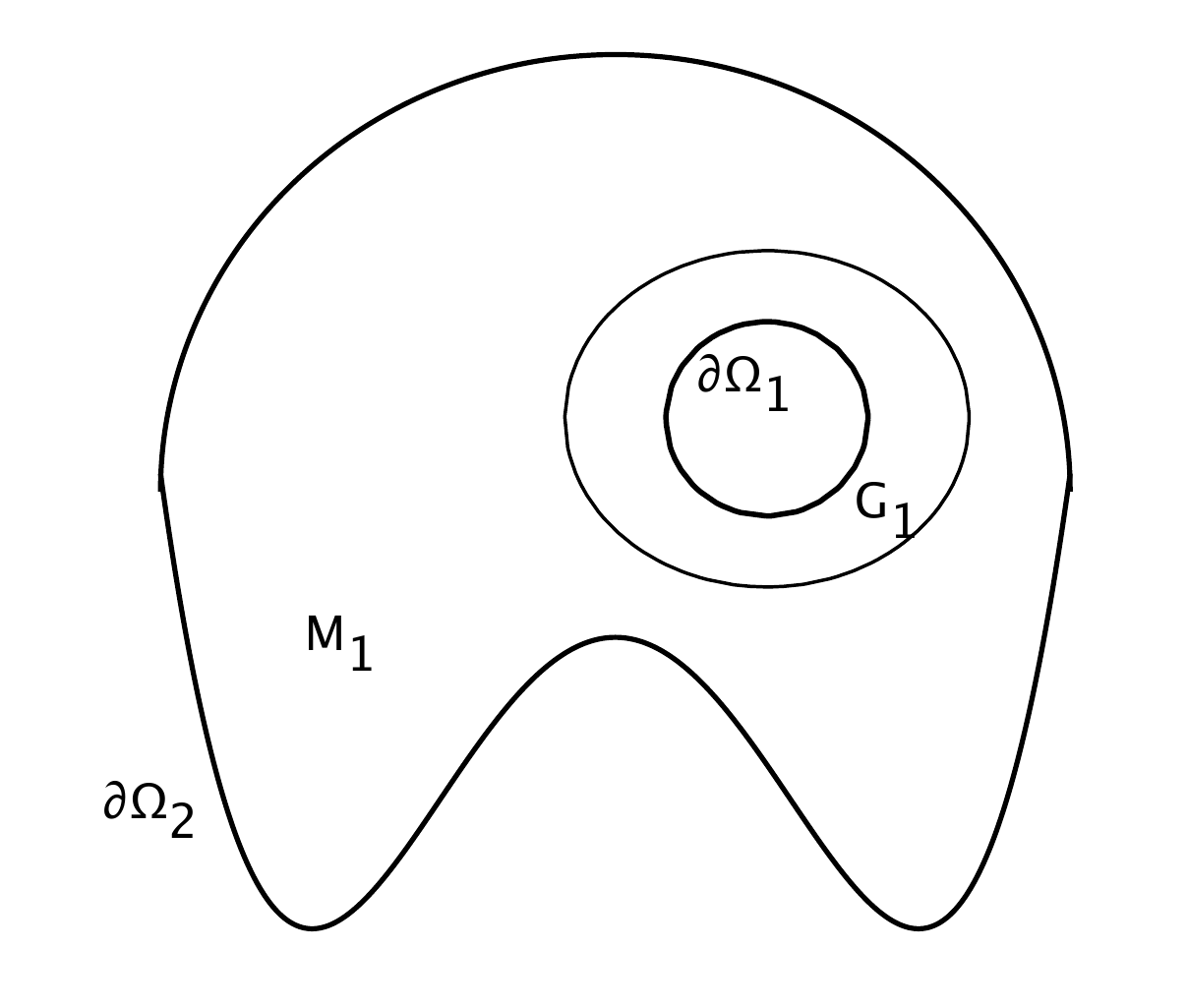}\end{center}   
     
     Let us note that 
  $K_i$  satisfies that, for all $x\in K_i$,
$\psi(x)=\Gamma_i \varphi(x)$ and $\grad \psi(x)=0$ and $\grad \varphi(x)=0$.
Indeed 
the first equality is true by continuity and the  other two because otherwise  $\psi=\Gamma_i \varphi$ in a neighborhood of that point, which contradicts the notion of boundary.
     
Moreover $\partial M_1=  \partial \Omega_2\cup K_1$ (and $\partial M_2=  \partial \Omega_1\cup K_2$).

      Indeed it is clear that $\partial M_1 \subset \partial \Omega_2\cup K_1$.  To prove the reverse inclusion, let $x\in K_1$, then  for all $r>0$,       $B(x,r) \cap (\R^N\setminus G_1) \neq \emptyset$.  Taking $r$ such that $B(x, r)\subset \Omega$ one gets  $x\in \overline{M_1}$. On the other hand $x\notin M_1$, since  if $x\in M_1$ there exists some ball $B(x, \epsilon) $ included in $M_1$, and $x\in K_1$ implies $B(x, \epsilon) \cap G_1\neq \emptyset$, a contradiction. 
      
     It is clear that $\partial \Omega_2$ is included in $\overline{M_1}$, and it  has no point of $M_1$ since $M_1$ is an open subset of $\Omega$. 
      
Let us admit for a while these three  simple claims,           
         
       {\bf Claim 1}       $M_i$ is connected for $i=1,2$.

      {\bf Claim 2} 
     $\partial M_i$ has at most two connected components i.e.       
$K_i$ is connected. 
    
    {\bf Claim 3}
       $\partial( M_1 \cap M_2) = K_1 \cup K_2$,

\noindent  and let us finish the proof.       We use first  claim 3 .   Since $\psi\in W^{2,p}$ for all $p< \infty$, using    Sard's theorem   in Sobolev spaces \cite{Fi},  there exists some constant $c_i$ such that $\psi_{\mid  K_i} = c_i$.

              Suppose that there exists one point  $\bar x\in M = M_1 \cap M_2$   such that $\psi (\bar x) < \min (c_1, c_2)$. Then $\psi $ would  have a local minimum inside $M$,   a contradiction with Hopf principle, see Remark \ref{remhopf}. Then  the minimum is achieved on the boundary of $M$, suppose   to fix the ideas that $c_1\leq c_2$.  Now take a ball in $M$  where $\psi>c_1$ that touches $K_1$ at some point $x_1$, by Hopf's principle $\grad \psi(x_1)\neq 0$ which contradicts the fact that, for all $x\in K_i$, $\grad \psi(x)=0$.

Now since $\psi$ cannot be locally constant we have proved that $M=\emptyset $ and then $\Gamma_1 = \Gamma_2$, hence $\Gamma = \gamma$.
This ends the proof of Theorem \ref{uni}.

\noindent Proof of  Claim 1 
              
              To fix the ideas we consider the case of $i=1$. 
If $M_1$ is not connected  there exists  at least two connected components of $M_1$; one denoted $M_{2,1}$ whose boundary contains $\partial \Omega_2$ and $\partial G_1$ and another connected component which we denote $M_{0,1}$.  We shall see that $M_{0,1}$ is simply connected. If not there exists some open regular domain $O^\prime$  with $\partial O^\prime \subset M_{0, 1}$ and $O^\prime$ is not included in $M_{0, 1}$. Then let $x\in O^\prime\setminus M_{0, 1}$, then  $x\in G_1$ and, since $G_1$ is connected,  $G_1\subset O^\prime$.
               
                We have obtained that $\partial G_1 \subset O^\prime$,  but this contradicts the fact that $\partial G_1\cap M_{2,1}\neq\emptyset $.

             Then   $M_{0,1}$ is an open set such that 
                $\partial M_{0,1}$ is connected and included in $K_1$. By Sard's theorem there exists a constant $c_{0,1}$ such that $\psi  = c_{0,1}$ on $\partial M_{0,1}$. Then one gets a contradiction with Hopf's principle (because either there exists a minimum inside $M_{0,1}$ and   this contradicts Remark \ref{remhopf}, or $c_{0,1}$ is a minimum for $\psi$ and taking some point on the 
       boundary  of $M_{0,1}$ which possesses the interior sphere condition, one gets once more a contradiction with Hopf).
                
                 We have obtained that $M_{0,1}= \emptyset$  and $M_1$ is connected. 
                 
\noindent Proof of Claim 2
                  
                  Suppose that $\partial M_1$ has at least three connected components. Then, since $\partial\Omega$ has two connected components, there exists some  domain  $O^\prime$  included in $\Omega$, such that $\partial O^\prime \subset M_1$ and $O^\prime $ is not included in $M_1$. Let $\bar x\in O^\prime \setminus M_1$.   Since  $\Omega \setminus \overline{M_1}\subset G_1$ and $G_1$ is connected,  then $G_1\subset O^\prime $.  This implies that $\partial \Omega_1 $ is also in the interior of $O^\prime$. Then  one cannot have $\Omega_1 \subset \subset \Omega_2$, a contradiction.

\noindent     Proof of Claim 3 
     
      We already have 
             \begin{eqnarray*}
             \partial( M_1 \cap M_2) &=& \overline{M_1\cap M_2} \setminus (M_1\cap M_2) \\
             &\subset&\overline{M_1} \cap \overline{M_2}\cap \left((\Omega\setminus M_1) \cup (\Omega\setminus M_2)\right)\\
             &\subset & ( \partial M_1 \cap \overline{M_2}) \cup (\partial M_2 \cap \overline{M_1})\\
            &\subset&  K_1 \cup K_2.
            \end{eqnarray*}   
To prove the reverse inclusion,      we  recall that, by hypothesis, $G_1\cap G_2 = \emptyset$. 
              
Then  $K_1\subset G_1 \subset M_2$. Let  $x\in K_1$ and $r>0$.  We need to prove that $B_r(x) \cap (M_1 \cap M_2) \neq \emptyset$.  Since $M_2$ is open there exists a ball $B_\epsilon (x) \subset M_2$. One can assume that $\epsilon < r$.  Since $x\in K_1 \subset  \partial M_1$, $B_\epsilon (x) \cap \overline{M_1} \neq \emptyset $ and also  $B_\epsilon (x) \cap M_1 \neq \emptyset $.
                Finally $B_\epsilon (x) \cap M_1 \cap M_2 \subset B_r(x) \cap M_1 \cap M_2$ hence $x\in \overline{M_1\cap M_2}$.  $x\in K_1$ implies $x\notin M_1$ hence $x\notin M_1\cap M_2$.  In the same manner $K_2 \subset \partial (M_1\cap M_2)$.

\bigskip
    We have also obtained the following strong comparison principle 
                  
                  \begin{theo}
   Suppose that  $N = 2$ and  that $\Omega$ and $\partial \Omega$ are connected. 
                   Suppose that $u$ and $v$ are  $ W^{2,p}(\Omega)$  for  some $p> 2$ respectively  super and sub solution of
                   $$\begin{array}{l} F(x, \nabla u, D^2 u)= f\ \mbox{in}\ \Omega\\
                   u=0\ \ \mbox{on}\ \partial\Omega
                   \end{array}$$
with $f\leq 0$  suppose that $u \geq v$ in $\Omega$. 
 
 Then either $u> v$ inside $\Omega$ or $u\equiv v$. 
\end{theo}

{\em Proof} : 
We     begin to prove that there exists a neighborhood of the boundary where either $u \equiv v$ or $u > v$ and in that last case one can conclude by using the  comparison principle in 
Theorem \ref{complambda},  and $u> v$ everywhere. 
We then consider the case where $u\equiv v$ on a neighborhood of $\partial\Omega$. As in the previous proof, let
$$A\mbox{ be the  connected component of  }\{ x, u(x) = v(x),  \ \nabla u (x) \neq 0, \ {\rm or} \ \nabla v(x) \neq 0\}$$
whose boundary contains $\partial\Omega$.
Let $G:=\overline{A}$,  $K := \partial G\cap\Omega$, $ M = \Omega \setminus G$.

We want to prove that $M=\emptyset$.
Suppose not, proceeding as above, it can be proved  that $M$  and $K$ are  connected, hence $\partial M= K$. Furthermore, on $K$, $\grad u=0$ and $\grad v=0$. 

Then, using Sard's theorem, 
 $ u = cte $ on $K$ and,   since $M$ cannot  contain  a local minimum   of $u$, the minimum of $u$ is achieved on the boundary i.e. on $K$; this contradicts the Hopf principle   on some point where the interior sphere condition is satisfied, since $\nabla u = 0$ on $K$.  Finally $M = \emptyset$ and $u = v$ in $\Omega$.

\section{Further results about the principal eigenvalues}

In all this section we suppose that  $F$ is given by (\ref{deff}) and $\tilde F$ satisfies (F) and (J) and $h$ satisfies  (H).

   \subsection{Properties  concerning  the dependance on the domain of the eigenvalues.}

From the definition of $\lambda^+$ and $\lambda^-$ it is clear that $\Omega_1\subset\Omega_2$ implies  $\lambda^+(\Omega_1)\geq \lambda^+(\Omega_2)$.
The next two theorems concern the strict  monotonicity of the principal eigenvalues with respect to the domain inclusion. We state them  in the case of the eigenvalue $\lambda^+ (\Omega),$ the symmetric results hold for ${\lambda^-} (\Omega)$ with obvious changes. In the sequel when no ambiguities arise we shall write $\lambda^+$ without writing the dependence on the domain.
\begin{theo}\label{zerob}
 Suppose that  $\Omega$ is a smooth bounded domain in $\R^N$. Suppose that $c(x)+\lambda^+> 0$ in $\Omega$. If $u$ is a positive solution of

$$F(x, \nabla u, D^2
u)+(c(x)+{\lambda}^+
)u^{1+\alpha} =0\, \ \rm{in}\  \Omega,$$ then there exists $(\partial\Omega)^\prime$ a connected component of $\partial\Omega$ such that
$$u=0\quad \mbox{on}\quad
(\partial\Omega)^\prime.$$
\end{theo}
{\em Proof.}
To begin with, let us note that $u$ must be zero somewhere  on $\partial \Omega$. Suppose not, i.e. $u>0$ on $\overline\Omega$. Using  the hypothesis $c(x)+\lambda^+ >0$, this implies  that there exists $\varepsilon>0$ and $\lambda^\prime> \lambda^+$ such that $u_\varepsilon:=u-\varepsilon>0$ is a solution of
 $$F(x, \nabla u_\varepsilon, D^2  u_\varepsilon)+( c(x)+\lambda^\prime ) u_\varepsilon^{1+\alpha}\leq 0,$$
 which   contradicts the definition of $\lambda^+$.

Let $v$ be an eigenfunction associated to $\lambda^+$,  so that $v>0$ in $\Omega$ and $v=0$ on  $\partial\Omega$.
 We now consider $\tau = \sup {v\over u}$ which exists and is finite. 
 
 Indeed, it is sufficient to prove it near the boundary. As seen in the proof of Theorem \ref{uni}
 there exists some constant $c_1$ such that $v\leq c_1 d$; on the other hand by Hopf's principle  there exists $c_o$ such that $u\geq  c_o d$; this implies that  $\tau$ is finite.

 Again reasoning as in  Theorem \ref{uni},  the infimum must be achieved at least on one point of  the boundary  in the sense that there
exists a sequence $x_n\in \Omega$, $x_n\rightarrow \bar
x\in\partial\Omega$ such that  ${v(x_n)\over u(x_n)}\rightarrow
\tau$.

Let $(\partial\Omega)^\prime$ be the connected component of $\partial\Omega$ that contains $\bar x$. Let $\Gamma^+ = \{ x\in (\partial\Omega)^\prime,u(x)>0\}$.
We want to prove that $\Gamma^+=\emptyset$.

 For that aim we consider  the set
 $$A = \{ x\in (\partial\Omega)^\prime,\ \displaystyle \limsup_{y\rightarrow x} {v(y)\over u(y)}  < \tau\}.$$
 We shall prove that $A$ is both closed and open in $(\partial \Omega)^\prime$.
 If we suppose that $\Gamma^+$ is non empty, it implies that $A$ is not empty since $ \Gamma^+\subset A$ and therefore $A=(\partial\Omega)^\prime$. 
We have obtained that $\displaystyle \limsup_{y\rightarrow x} {v (y) \over u (y) }< \tau$ everywhere on $(\partial\Omega)^\prime$, contradicting the previous observation.

By Hopf principle and the regularity result one has  $|\nabla v|\geq L_1>0$,
 on a neighborhood of the boundary. 
In the sequel  for $x\in \partial \Omega$, we denote as $B_{x}(r_1)$ some ball such that $B_{x} (r_1) \subset \Omega$ and  $\overline{B_{x} (r_1)} \cap \partial{\Omega} =x $.  The existence of such ball is implied in particular by the assumption that $\partial \Omega$ is ${\cal C}^2$.  ($W^{2, \infty}$ is sufficient).

\noindent{\bf Claim}  : The complementary of $A$ is closed in $(\partial \Omega)^\prime$.

We want to prove that if $x_n$ is some sequence in $(\partial\Omega)^\prime$
which converges to $x$ such that $\displaystyle \limsup_{y\rightarrow x_n}{v(y) \over u(y) }=\tau$, then

$$    \limsup_{y\rightarrow x} {v(y) \over u(y) }  =\tau.$$
By the strict comparison principle Corollary \ref{234} one has that there exists $r_1$ independent on $n$ such that  $v = \tau u$ in $B_{x_n} (r_1)$. 

Since  $B_{x_n}(r_1)\cap B_x(r_1)\ \neq
\emptyset$ for $n$ large enough, there exists $y_n\in B_{x}(r_1)$,
$v(y_n)=  \tau u(y_n)$,  again using Corollary \ref{234} we get that $v= \tau u$ in all the ball $B_{x}(r_1)$  and then
$x$ does not belong to $A$.

\noindent{\bf Claim }: The set $A$  is closed.

Let $x_n\in \partial \Omega $ such that for all $n$,  $\displaystyle\limsup_{y\rightarrow x_n} {v (y)\over u (y)}< \tau$
  and $x_n$
converges to $\bar x\in (\partial\Omega)^\prime$, let us prove that $\bar x\in A$. Using Proposition \ref{prop1},  there exists $r_1$ such that 
$v<   \tau u$ in $B_{x_n}(r_1)$ , and since $B_{x_n}
(r_1)\cap B_{\bar x} (r_1)\cap \Omega\neq \emptyset, $ there exists $y\in B_{\bar x}(r_1)$ such that $v (y)<  \tau u(y)$,   which  implies
that
$v< \tau u$ in all the ball $B_{\bar x}(r_1)$.   

There are two cases either $v(\bar x) = 0 < u(\bar x)$, then $\bar x\in A$ and we are done. Or   $v(\bar x) = \tau u(\bar x)=0$ and
we need to prove that 
$\displaystyle\limsup_{y\rightarrow \bar x } {v(y) \over u(y) }  <\tau.$
If not by Corollary \ref{234}, $v = \tau  u$ in the whole ball  $B_{\bar x}(r_1)$, a contradiction, so $\bar x\in A$.
This ends the proof. 

\bigskip
It is well known (and easy to prove) that if $\Omega^\prime\subset\Omega$ then ${\lambda}^\pm (\Omega^\prime)\geq {\lambda}^\pm (\Omega).$ In the next Theorem we shall prove that the monotonicity is strict.

\begin{theo} \label{thf} Let $\Omega$ be a domain with a connected boundary, and $\Omega^\prime$ a subdomain of $\Omega$. 

Suppose that $\partial \Omega$ is not included in $\partial\Omega^\prime$ or, for $N=2$, suppose only that $\overline{\Omega^\prime}\neq\overline\Omega$. Then
$${\lambda}^\pm (\Omega^\prime)> {\lambda}^\pm (\Omega).$$
\end{theo}
{\em Proof:} Let us begin by the first case i.e. we suppose that there exists 
 $x_0\in \partial \Omega$ such that $x_0\not\in \partial\Omega^\prime$ and hence, there exists 
 $\delta>0$, such that
$B(x_0, \delta)\subset \R^N\setminus \overline{\Omega^\prime}$.

Let $\Omega^{\prime\prime}$ be a smooth domain whose boundary has only one connected component and such that 
$$\Omega^\prime\subset\Omega^{\prime\prime}\subset\Omega \setminus
\overline{ B (x_0,\delta)}.$$
We know that 

$$\lambda^+(\Omega^\prime)\geq\lambda^+(\Omega^{\prime\prime})\geq\lambda^+(\Omega \setminus
\overline{ B (x_0,\delta)})\geq \lambda^+(\Omega).$$

Suppose by contradiction that ${\lambda}^+ (\Omega^\prime)={\lambda}^+ (\Omega)$; this implies that $\lambda^+(\Omega^{\prime\prime})={\lambda}^+ (\Omega):=\lambda^+$.

 Let  $v$ be  an eigenfunction for
$\Omega$. So in particular it is positive in $\Omega$. On the other hand it is also a solution
of
$$F(x, \nabla \varphi, D^2
\varphi)+(c(x)+{\lambda}^+  )\varphi^{1+\alpha} =0,
$$
in $\Omega^{\prime\prime}$,  then using  Theorem \ref{zerob}, this  would imply that $v=0$ on $\partial
\Omega^{\prime\prime}$, a contradiction with the fact that $v>0$ in $\partial \Omega^{\prime\prime} \cap \Omega$. This ends the first case.

We are left to prove the case $N=2$ and $\partial\Omega\cap\partial\Omega^\prime=\partial\Omega$. 
By contradiction we shall suppose that $\lambda^+(\Omega^\prime)=\lambda^+(\Omega)$. We can assume that $\Omega^\prime$ is regular by replacing $\Omega^\prime$ by some  regular subset  $\Omega^{\prime\prime}$ such that 
$\partial \Omega^{\prime\prime}  \cap  \partial\Omega=\partial\Omega$ which contains $\Omega^\prime$ since it satisfies also $\lambda^+(\Omega^{\prime\prime}) = \lambda(\Omega)$. For simplicity we rename this set    $\Omega^\prime$  and we consider the respective eigenfunctions $\phi$ and 
$\phi^{\prime}$.  By hypothesis they satisfy the same equation in $\Omega^\prime$.

Let 
$$\tau=\sup_{\Omega^\prime}\frac{\phi^\prime}{\phi},$$ 
and, proceeding as in the proof of Theorem \ref{uni}, it is possible to prove that $\tau$ is bounded and there exists $G$ a closed connected neighborhood of $\partial\Omega$ where $\phi^\prime=\tau\phi$ and such that 
$$\grad\phi=\grad\phi^\prime=0\mbox{ in }\ K:=\partial G\cap\Omega.$$
Hence using Sard's theorem for functions in $W^{2,p}$, there exists a constant $c$ such that

$$\phi=c,\ \phi^\prime=\tau c, \mbox{ in}\  K.$$ 
This of course leads to a contradiction, because either $\phi<c$ somewhere  in $\Omega\setminus G$ and then $\phi$ would have a local minimum which contradicts Remark \ref{remhopf}. Or $\phi\geq c$ in $\Omega\setminus G$ and then by Hopf's Lemma and Lemma \ref{sph} there is a point where $\grad \phi\neq 0$ in $K$ which is again a contradiction.
This ends the proof.
\bigskip

In the next theorem we suppose, to fix the ideas ${\lambda}^+ <{\lambda}^-$, with obvious symmetric results in the other case. 
 
 \begin{theo}\label{lambda-}
 
  Suppose that  $\Omega$ is  a bounded  regular domain. Suppose that $f\geq 0$, then there exists no solution $u$ for the equation    
  $$\left\{ \begin{array}{lc}
  F(x, \nabla u, D^2 u)+(c(x)+{\lambda}^+  )u^{1+\alpha} = f& {\rm in} \ \Omega\\
 u=  0& \ {\rm  on} \ \partial \Omega.
 \end{array}\right.$$
such that $u(x_o)>0$ for some $x_o\in\Omega$.

  Suppose in addition that $\partial \Omega$ is connected, and that  $f<0$ somewhere near the boundary, then there is no solution  $u$  
      $$\left\{ \begin{array}{lc}
  F(x, \nabla u, D^2 u)+(c(x)+{\lambda}^-  )|u|^\alpha u  = f& {\rm in} \ \Omega\\
 u=  0& \ {\rm  on} \ \partial \Omega
 \end{array}\right.$$ 
such that $u(x_o)<0$ for some $x_o\in\Omega$.
  \end{theo}
{\em Proof of Theorem \ref{lambda-}}:
The first part is a mere application of the minimum principle : Since $f\geq 0$ and $ {\lambda}^+< {\lambda}^-$, and  since $ u = 0$ on the boundary, $ u\leq 0$, which contradicts the assumption. 
  
    We now prove the second part. 
    
    Since $u(x_o)<0$ for some $x_o$ let  $\varphi^-$ be some normalized eigenfunction  and 
   $\Gamma = \inf\{ t, \ |u|\leq t|\varphi^-|\}$. 
   Then $\Gamma >0$ and 
   $|u|\leq \Gamma |\varphi^-|$.    
     We prove that the supremum $\Gamma = \sup {u\over \varphi^-}$ is achieved on the boundary. If not there exists some compact set 
    $K$ large enough in order that ${u\over \varphi^-}\leq \Gamma^\prime$ on $\Omega \setminus K$ and $\Gamma^\prime < \Gamma$.

     Let $K$ be some compact set such that  ${\lambda}^+(\Omega\setminus K) > {\lambda}^+ (\Omega)$.  Then using the comparison principle on $\Omega \setminus K$, since $u$ is a supersolution  and  $u\geq \Gamma^\prime \varphi$ on $\partial K$,and then 
     $u \geq \Gamma^\prime \varphi^-$ in $K$, finally $u \geq \Gamma^\prime \varphi^-$ in the whole of $\Omega$,  which yields a contradiction. 
     
      Then the supremum is "achieved" on the boundary and then  the strict comparison principle  in Proposition  \ref{prop1} implies that $u = \Gamma \varphi$ around the boundary. In particular one gets $f \equiv 0$ around the boundary, which is once more a contradiction.

  \subsection{Further properties}
We want to prove that we can recover some of the standard properties of eigenvalues for linear elliptic equations. We  consider the Dirichlet problem 
 \begin{equation} \label{oev}
 \left\{ \begin{array}{lc}
  F(x, \nabla u, D^2 u)+(c(x)+{\lambda}  )u^{1+\alpha} = 0& {\rm in} \ \Omega\\
 u=  0& \ {\rm  on} \ \partial \Omega.
 \end{array}\right.
 \end{equation}
Non trivial solutions of (\ref{oev}) will be called eigenfunctions.

We now recall the following result, which is an easy consequence of the definitions of $\lambda^{\pm}$ and the maximum and the minimum principle (see \cite{BD5},\cite{BD6})
 \begin{theo}\label{th6} Suppose that $\Omega$ is a bounded regular domain,

 1) If $\lambda > \lambda_1=\sup ({\lambda}^+(\Omega),{\lambda}^-(\Omega))$, then every non trivial solution of (\ref{oev}) changes sign in $\Omega$.

 2) For any $\lambda$ between ${\lambda}^+$ and ${\lambda}^-$ there are no nontrivial solutions of  (\ref{oev})

  \end{theo}
As a consequence of the previous result, we now have further results regarding the signs of the eigenfunctions.
 \begin{theo}\label{th7}  Suppose that  $\Omega$ is some  bounded smooth domain, 
 and suppose that one of the following holds
 
 \begin{enumerate}
 \item $\lambda^+=\lambda^-$ and $\partial\Omega$ is connected
 \item $N=2$, $\lambda^+=\lambda^-$, $\partial\Omega$ has at most two connected components
 \item $N=2$, $\partial\Omega$ is connected
 \end{enumerate}
 then any eigenfunction  corresponding to  $\lambda=\lambda^\pm$ is of constant sign.
 \end{theo}
 An application of the previous result is the following
 \begin{cor}\label{cor46}
 In the hypothesis of Theorem \ref{th7} the eigenvalues $\lambda^+$ or $\lambda^-$ are isolated i.e.  there exists $\delta>0$ such that for  any $\lambda\in ]\lambda^\pm,\lambda^\pm +\delta[$, the solutions of (\ref{oev}) are trivial.
\end{cor}
 {\em Proof of Theorem \ref{th7}.}
 Let $\lambda_1 = \lambda^+ = \lambda^-$ and suppose that $\partial\Omega$ is connected.
    Suppose by contradiction that there exists a solution of (\ref{oev}) which changes sign. We define
 $\Omega^+=\{x\in\Omega;\ u(x)>0\}$ and  $\Omega^-=\{x\in\Omega;\ u(x)<0\}$. Then clearly, for any $\tilde\Omega^+$  (respectively $\tilde\Omega^-$) connected component of $\Omega^+$ (respectively $\Omega^-$):
$$\lambda^+(\tilde\Omega^+)=\lambda^-(\tilde\Omega^-)=\lambda_1.$$

$\partial\tilde\Omega^+\cap\partial\Omega\neq\partial\Omega$ is not possible  since it would imply, by Theorem \ref{thf} $\lambda^+(\tilde\Omega^+)>\lambda^+(\Omega)$. But, on the other hand, if
 $\partial\tilde\Omega^+\cap\partial\Omega=\partial\Omega$, then  $\partial\tilde\Omega^-\cap\partial\Omega\neq\partial\Omega$, and then the contradiction is given by the fact that it would imply
  $\lambda^-(\tilde\Omega^-)>\lambda^-(\Omega)$.

 We now suppose to be in the second case, in particular we suppose that  $\partial \Omega$ has  two connected components. For $i=1,2$, $\Omega_i$ denote two open simply connected sets such that 
 $\Omega_1 \subset \subset \Omega_2$ and $\Omega = \Omega_2\setminus \overline{\Omega_1}$.
 
 Let $\varphi$ be a positive eigenfunction in $\Omega$ and $v$  an eigenfunction which changes sign. Let  $\Omega^\pm$ be the set where $v$ is positive (respectively negative). If $\partial \Omega^+ \cap \partial \Omega_i\neq \emptyset$ then reasoning as in the proof of Theorem \ref{zerob}, $\partial \Omega^+ \cap \partial \Omega = \partial \Omega_i$. 
 The same is true for $\Omega^-$, hence one can assume without loss of generality that $\partial \Omega^+ \cap \partial \Omega = \partial \Omega_1$ and  $\partial \Omega^-\cap \partial \Omega = \partial \Omega_2$, finally $\partial \Omega ^+ \cap \Omega = \partial \Omega ^-\cap \Omega\neq \emptyset$. 
 
 Since $\lambda^+=\lambda^-$, $-\phi$ is an eigenfunction corresponding to $\lambda^-$. Reasoning as in the proof of Theorem \ref{uni} we can define $\Gamma =\sup {v\over \varphi}$ which exists according to the estimates on $v$ near the boundary and it is "achieved" on $\partial \Omega_1$, and $\gamma = \sup {-v\over \varphi}$ "achieved " on $\partial \Omega_2$, and prove that 
 $v = \Gamma \varphi$ on a neighborhood of $\partial \Omega_1$ and 
$v = -\gamma \varphi$ on a neighborhood of $\partial \Omega_2$. 

One define as in the proof of Theorem \ref{uni}  $M_1$ , $K_1$, $M_2$, $K_2$ and one has on $ K_1$  and on $K_2$, $\nabla v = \nabla \varphi = 0$. Since $\varphi >0$ in $M_1 \cap M_2$, using Sard's theorem $c_i = \varphi_{\mid_{K_i}}$.  Again reasoning as in Theorem \ref{uni},
$ \varphi$ must achieve  its infimum on $K_1$ or $K_2$ where its gradient is zero, a contradiction with Hopf principle. 
 \bigskip

We consider the last  case, i.e. $N= 2$, $\lambda^+ \neq \lambda^-$ and $\partial \Omega$ is connected.    It is clear  with the maximum principle that the result is  true for the smaller of the two eigenvalues. We suppose to fix the ideas that $\lambda^-< \lambda^+$ and we prove that every eigenfunction corresponding to $\lambda^+$ is positive.   
   
    We denote by $\varphi$ some positive eigenfunction for $\lambda^+(\Omega)$ and by $\psi$ some eigenfunction which changes sign. 
     We denote by $\Omega^+$ and $\Omega^-$ the sets   where $\psi$ is respectively positive and negative.  
This implies that $\lambda^+ (\Omega^+) = \lambda^+(\Omega)$ and $\lambda^-(\Omega^-)= \lambda^-(\Omega)$; but, in dimension $N=2$, this contradicts Theorem \ref{th7}.

 \bigskip
 {\em Proof of Corollary \ref{cor46}.} The result needs to be proved for $\lambda_1=\sup(\lambda^+,\lambda^-)$.
  Suppose  by contradiction that there exists a sequence of eigenvalues $\lambda_n$, such that $\lambda_n\rightarrow\lambda_1$, $\lambda_n>\lambda_1$.

  Let $u_n$ be a sequence of  solutions of (\ref{oev}) with $\lambda=\lambda_n$ such that $|u_n|_\infty=1$.  This implies that the H\"older's  norm is uniformly bounded with respect to $n$ (see Proposition \ref{prophold}).

  Then  $(u_n)_n$ is relatively compact and it converges up to a subsequence in ${\cal C} (\bar\Omega)$ towards a solution $u$ of (\ref{oev}) with $\lambda=\lambda_1$.

  By Theorem \ref{th7}, $u$ must be either positive or negative, which implies that for $n$ large enough,  for any $K=\overline\Omega_1\subset\Omega$, for $n$ large enough, $u_n$ has constant  sign in $K$. 

Without loss of generality we can suppose that $\lambda_1=\lambda^+$ and $u>0$ in $\Omega$ and hence $u_n$ is positive in $K$.
 
 We choose $\Omega_1$  a regular  subset in $ \Omega $ such that 
 $\lambda^-(\Omega\setminus\overline{\Omega_1})>\lambda_n$.

By minimum principle, since 
$$u_n\geq 0,\ \mbox{in }\ \partial( \Omega\setminus\overline{\Omega_1})$$
it implies that $u_n\geq 0$ in $\Omega\setminus\overline{\Omega_1}$ and then $u_n$ is positive in $\Omega$ which contradicts Theorem \ref{th6}.

      \end{document}